\documentclass[a4paper]{article}
\usepackage[comma,numbers,square,sort&compress]{natbib}
\newtheorem{lemma}{Lemma}
\newtheorem{theorem}{Theorem}

\newtheorem{remark}{Remark}
\newtheorem{assumption}{Assumption}

\usepackage{graphicx,amsmath,amssymb}
\usepackage{subfigure,pgf,tikz,pst-node}
\usetikzlibrary{arrows,shapes,automata,positioning,fit}

\begin{document}
\title{Robust Consensus Tracking of Heterogeneous Multi-Agent Systems under Switching Topologies}
\author{Yutao Tang \footnote{Y. Tang is with School of Automation, Beijing University of Posts and Telecommunications, Beijing 100876, China.  E-mail\,$: yttang@amss.ac.cn$}}
\date{}
\maketitle

{\noindent\bf Abstract} {In this paper, we consider a robust consensus tracking problem of heterogeneous multi-agent systems with time-varying interconnection topologies. Based on common Lyapunov function and internal model techniques, both state and output feedback control laws are derived to solve this problem.  The proposed design is robust by admitting some parameter uncertainties in the multi-agent system.}


\section{Introduction}

The past decade has witnessed a rapid development in the field of multi-agent system and many fruitful results have been obtained. As one important topic,  the consensus problem is actively studied due to its
numerous applications such as cooperative control of unmanned aerial vehicles, communication among sensor networks, and formation of mobile robots (see \cite{fax2004information, ren2008distributed} and the references therein). Consensus means that a group of agents reaches an agreement on a physical quantity of interest by interacting with their local neighbors. Roughly speaking, existing consensus problems can be categorized into two types: consensus without a leader (i.e., leaderless consensus or synchronization) and consensus with a leader. The case of consensus with a leader is also called leader-following consensus or consensus tracking.

The consensus tracking problem of multi-agent systems has been studied for many years. This problem with agents in the form of single-integrators or double-integrators has been widely investigated (\cite{ren2008distributed, meng2011leaderless, li2013quantized}).  In \cite{hong2006tracking}, the authors proposed a distributed observer-based control law using local information to track an integrator-type leader. Later, this work has been extended to multi-agent systems with a general linear dynamics (\cite{hong2009multi, ni2010leader}) under switching topologies. Other extensions including consensus tracking with an unknown-input driven leader or a designed abstracted system have been studied in \cite{tang2014leader,tang2013hierarchical,tang2015auto}

Recently, a general framework based on output regulation theory (\cite{wonham1979}) has been developed for multi-agent consensus (\cite{wang2010distributed, su2013general}), where distributed control could achieve both asymptotic tracking and disturbance rejection. For example, as an extension to the results in \cite{hong2009multi},  a cooperative output regulation problem of heterogeneous multi-agent systems was formulated in \cite{su2012tac} and solved by devising a distributed observer, where the exosystem can be viewed as the leader in the leader-following formulation. However, these results were obtained with a fixed topology. A similar problem but under a class of switching topologies was considered in \cite{su2012smc-cooperative}. Both dynamic state feedback and measurement output feedback control laws were proposed by using the full information of the exosystem which may be undesirable sometimes. Without requiring exactly known system matrices, many robust results were further obtained which can handle parameter uncertainties in the multi-agent system.  For example, the authors in \cite{wang2010distributed} proposed a robust control law based on internal model design with a no-loop graph which admits parameter uncertainties in the system matrices. Later, the no-loop graph assumption was removed in \cite{su2013general}.  These results were mainly achieved with fixed topology. Very few results were obtained on robust output regulation with switching topologies, except that a control design was investigated in \cite{wang2011robust} based on canonical internal models to deal with small parameter uncertainties. However, this method was only applicable to the square follower system and led to an output feedback control law.

Our main motivation is to investigate the robust consensus tracking problem of heterogeneous multi-agent systems under switching topologies, which may have parameter uncertainties in their system matrices. To facilitate the control goal, a distributed observer regarding switching topologies is at first constructed for each agent, and thus the robust consensus tracking problem can be solved by incorporating an internal model in the controller for each agent. Our main contributions are at least two-fold:
\begin{itemize}
  \item We extend the conventional leader-following consensus (e.g., \cite{hong2006tracking, ren2008distributed}) to general linear multi-agent systems with parameter uncertainties. When there are no such uncertainties, these results are consistent with existing consensus and distributed output regulation results, e.g., \cite{hong2006tracking, su2012tac}.
  \item We consider the distributed/cooperative robust output regulation problem under switching topologies, and robust control laws are proposed for heterogeneous multi-agent systems, while many existing regulation results were derived for only homogeneous multi-agent systems (\cite{su2013general}) or for fixed graph cases (\cite{wang2010distributed, su2012tac}).
\end{itemize}
Moreover, both output and state feedback control laws are proposed to solve the the robust output regulation problem considered in \cite{wang2011robust}, and these control laws are also applicable to the follower system whose inputs have different dimensions with its outputs.

The rest of this paper is organized as follows. In Section 2, the problem formulation is given. Main results are presented in Section 3, where two types of robust consensus control laws are proposed. Finally, simulations and our concluding remarks are presented at the end.

Notations:  Let $\mathbb{R}^n$ be the $n$-dimensional Euclidean space, $\mathbb{R}^{n\times m}$ be the set of $n\times m$ real matrices, $\mathbf{0}_{n \times m}$ represents an $n\times m$ zero matrix. $\mathbf{1}_{n}$ represents a column vector of $n$ entries equal to 1. $\text{diag}\{b_1,{\dots},b_n\}$ denotes an $n\times n$ diagonal matrix with diagonal elements $b_i\; (i=1,{\dots},n)$; $\text{col}(a_1,{\dots},a_n) = [a_1^T,{\dots},a_n^T]^T$ for any column vectors $a_i\; (i=1,{\dots},n)$, $\mbox{vec}(X)=\mbox{col}(X_1,\dots,X_n)$ with $X_i$ as the $i$-th column of a $n$-column matrix $X$. A weighted directed graph (or weighted digraph) $\mathcal {G}=(\mathcal {N}, \mathcal {E}, \mathcal{A})$ is defined as follows, where $\mathcal{N}=\{1,{\dots},n\}$ is the set of nodes, $\mathcal {E}\subset \mathcal{N}\times \mathcal{N}$ is the set of edges, and $\mathcal{A}\in \mathbb{R}^{n\times n}$ is a weighted adjacency matrix (\cite{godsil2001}). $(i,j)\in \mathcal{E}$ denotes an edge leaving from node $i$ and entering node $j$. The weighted adjacency matrix of this digraph $\mathcal {G}$ is described by $A=[a_{ij}]_{i,\,j=1,\dots,n}$, where $a_{ii}=0$ and $a_{ij}\geq 0$ ($a_{ij}>0$ if and only if there is an edge from agent $j$ to agent $i$). A path in graph $\mathcal {G}$ is an alternating sequence $i_{1}e_{1}i_{2}e_{2}{\cdots}e_{k-1}i_{k}$ of nodes $i_{l}$ and edges $e_{m}=(i_{m},i_{m+1})\in\mathcal {E}$ for $l=1,2,{\dots},k$. If there exists a path from node $i$ to node $j$ then node $i$ is said to be reachable from node $j$. The neighbor set of agent $i$ is defined as $\mathcal{N}_i=\{j\colon (j,i)\in \mathcal {E} \}$ for $i=1,...,n$.  A graph is said to be undirected if $a_{ij}=a_{ji}$ ($i,j=1,{\dots},n$). The weighted Laplacian $L=[l_{ij}]\in \mathbb{R}^{n\times n}$ of graph $\mathcal{G}$ is defined as $l_{ii}=\sum_{j\neq i}a_{ij}$ and $l_{ij}=-a_{ij} (j\neq i)$.

\section{Problem Formulation}
Consider a group of $N+1$ agents, and $N$ of them are followers described by
\begin{align}\label{sys:follower}
\begin{split}
  \dot{x}_i&=(A_i+\Delta A_i)x_i+(B_i+\Delta B_i)u_i\\
  y_i&=(C_i+\Delta C_i) x_i+(D_i+\Delta D_i)u_i,\quad i=1,\dots,N
  \end{split}
\end{align}
where $x_i\in \mathbb{R}^{n_i},\, u_i\in\mathbb{R}^{p_i},\, y_i\in\mathbb{R}^{q}$ are the state, input, and output of agent $i$, respectively. $A_i\in \mathbb{R}^{n_{i}\times n_{i}}$,\,$B_i\in \mathbb{R}^{n_{i}\times {p_i}}$,\,$C_i\in \mathbb{R}^{q\times n_i}$,\,$D_i\in \mathbb{R}^{q\times p_i}$ represent the nominal part of the plant, while $\Delta A_i$,\,$\Delta B_i$,\,$\Delta C_i$,\,$\Delta D_i$ represent the uncertain part. Without loss of generality, assume the triple $(C_i,\,A_i,\,B_i)$ is controllable and observable.

Let $w_i=\mbox{vec}\left(\begin{bmatrix}
  \Delta A_i&\Delta B_i\\
  \Delta C_i&\Delta D_i
\end{bmatrix}\right)$ be the uncertain vector for agent $i$ and denote $A_i(w)=A_i+\Delta A_i$,\,$B_i(w)=B_i+\Delta B_i$,\,$C_i(w)=C_i+\Delta C_i$,\,$D_i(w)=D_i+\Delta D_i$,\,$w=\text{col}(w_1,\dots,w_N)$.

The leader(node $0$) is as
\begin{align}\label{sys:leader}
\begin{split}
\dot{v}=Sv,\quad y_0=-Fv
\end{split}
\end{align}
where $S\in \mathbb{R}^{m\times m}$,\,$F\in \mathbb{R}^{q\times m}$, and $(F, S)$ is observable.

Associated with this multi-agent system, a dynamic digraph $\mathcal{G}$ can be defined with the nodes $\mathcal{N}=\{0,1,..., N\}$ to describe the communication topology, which may be switching. If the control $u_i$ can get access to the information of agent $j$ at time instant $t$, there is an arc $(j,i)$ in the graph $\mathcal{G}$, i.e., $a_{ij}>0$.  Also note that $a_{0i}=0$ for $i=1,...,N$ since the leader won't receive any information from the followers. Denote the induced subgraph associated with all followers as $\bar{\mathcal{G}}$. The associated Laplacian of this digraph can be partitioned as
\begin{align*}
  L=\left[\begin{array}{c|c}
    0&\textbf{0}_{1\times N}\\ \hline
    \hat L&H
  \end{array}\right]
\end{align*}
where $\hat L\in \mathbb{R}^{N\times 1}$ and $H\in \mathbb{R}^{N\times N}$.

To achieve the coordination of the multi-agent system, node $0$ in $\mathcal G$ should be globally reachable, and therefore at least one agent in each component of ${\mathcal G}$ is connected to the leader. Otherwise, the coordination between the agents and the leader cannot be achieved. We say a communication graph is connected if the leader (node 0) is reachable from any other node of $\mathcal{G}$ and the induced subgraph $\bar{\mathcal{G}}$ is undirected. By Lemma 3 in \cite{hong2006tracking}, $H$ is positive definite if the communication graph is connected. Denote its eigenvalues as $\lambda_1\geq\lambda_2\geq{\cdots}\geq \lambda_N>0$.

In multi-agent systems, the connectivity graph $\mathcal{G}$ may be time-varying. To describe the variable interconnection topology, we denote all possible communication graphs as $\mathcal{G}_1$,\ldots,$\mathcal{G}_\kappa$, $\mathcal{P}=\{1,\dots, \kappa\}$, and define a switching signal $\sigma: [0,\infty)\rightarrow \mathcal{P}$, which is piece-wise constant defined on an infinite sequence of nonempty, bounded, and continuous time-intervals. Assume  $t_{i+1}-t_i\geq \tau_0, \;\forall i$, where $t_{i}$ is the $i$-th switching instant and $t_0=0$. Here $\tau_0$ is often called the dwell-time. Therefore, $\mathcal{N}_i$ and the connection weight $a_{ij}\; (i,j=0,1,\ldots,N)$ are time-varying. Moreover,  the Laplacian $L_{\sigma(t)}$ and the matrix $H_{\sigma(t)}$ associated with the switching interconnection graph $\mathcal{G}_{\sigma(t)}$ are also time-varying (switched at $t_i,\; i=0,1,\ldots$), though they are time-invariant in each interval $[t_{i},t_{i+1})$.

The following assumption on the communication graph is made.
\begin{assumption}\label{ass:graph}
The graph is switching among a group of connected graphs with the leader as its root.
\end{assumption}

The robust consensus tracking problem of this heterogeneous multi-agent system can be formulated as follows. {\it Given the multi-agent system composed of the follower \eqref{sys:follower}, the leader \eqref{sys:leader} and its corresponding communication graph $\mathcal{G}_{\sigma(t)}$, find a robust control law such that, there exists an open neighborhood $W$ of $w=0$, for any initial condition $x_i(0)\in \mathbb{R}^{n_i}$, $v(0)\in\mathbb{R}^{m}$, $w\in W$, the consensus tracking goal is achieved, i.e., $e_i=y_i-y_0 \to 0$ for $i=1,\dots,N$ as $t\to\infty$.}

Note that, the system of agents considered are all with general linear models and can be seen as extensions to the well-studied consensus tracking problem in existing results for integrator-type agents (\cite{fax2004information, ren2008distributed, hong2006tracking}).

Another assumption is needed to solve this robust consensus tracking problem.
\begin{assumption}\label{ass:trans}
Let $\lambda$ be any eigenvalue of $S$, it holds
      \begin{align}\label{eq:trans}
      \text{rank}\begin{bmatrix}
        A_i-\lambda I&    B_i\\
        C_i          &    D_i
      \end{bmatrix}=n_i+q.
      \end{align}
\end{assumption}
\begin{remark}
  This condition is often called the transmission zero condition (\cite{huang2004}) and is a key ingredient to solve the robust output regulation problem in centralized cases.
\end{remark}

\begin{remark}\label{rem:su-robust}
  When all agents have the same nominal system matrices, i.e., $A_i=A$,\,$B_i=B$,\,$C_i=C$,\,$D_i=D$, a similar result have obtained in \cite{su2013general} under a fixed topology. Note that the interconnection in our formulation allows switching topologies, which renders this problem more challenging.
\end{remark}

\section{Main Results}

In this section, we give a robust design scheme based on internal model techniques.

Unlike in decentralized cases (\cite{xi2007global}), we do not assume the availability of $y_0$ to all agents in our problem. An agent can get access to $y_0$ unless there is an arc from the leader to this agent. Since not all followers can directly get access to the leader's information, we first build a distributed observer for each agent as follows.
\begin{align}\label{sys:observer}
\dot{\eta}_i=S\eta_i-L_0F\eta_{vi}
\end{align}
where $\eta_{vi}=\sum_{j=0}^N a_{ij}(t)(\eta_i-\eta_j)$, $\eta_{0}=v$ and $L_0$ is a gain matrix to be selected.

To establish the performance of this distributed observer, the following lemma will be used.
\begin{lemma}\label{lem:lmi}
Let $P$ be a positive definite symmetric solution of the Lyapunov inequality
\begin{equation}\label{eq:lyapunov}
PS+S^TP-2F^TF<0.
\end{equation}
Take $L_0=\mu P^{-1}F^T$, then under Assumption \ref{ass:graph}, there exist positive constants $\mu^*$ and $c$, such that, when $\mu\geq \mu^*$, it holds
\begin{equation}\label{eq:lmi}
(S-\lambda_i^p L_0F)^T P+ P(S-\lambda_i^p L_0F)\leq -cP
\end{equation}
where $\lambda_i^p>0$ for $i=1,...,N$ are the eigenvalues of $H_p$ ($p\in \mathcal{P}$).
\end{lemma}
{\bf Proof}: Note that there are only finite graphs satisfying Assumption \ref{ass:graph}, the minimum eigenvalue of $H_p$ for all $p$ is well-defined and denoted as $\bar \lambda>0$. Taking $\mu^*=\max\{\frac{1}{\bar \lambda},1\}$, we obtain
\begin{align*}
  &(S-\lambda_i^p L_0F)^T P+ P(S-\lambda_i^p L_0F)\\
  &=S^TP+PS-2\mu\lambda_i^p F^TF\\
  &=S^TP+PS-2F^TF-2(\mu\lambda_i^p-1) F^TF\\
  &\leq S^TP+PS-2F^TF
\end{align*}
Since $S^TP+PS-2F^TF$ is negative definite, there exists a sufficiently small constant $c$ satisfying the inequality \eqref{eq:lmi}. \hfill\rule{4pt}{8pt}

The next lemma shows the convergence of \eqref{sys:observer} and thus guarantees that each agent can asymptotically get the state information of the leader.

\begin{lemma}\label{lem:observer}
  Under Assumption \ref{ass:graph}, there exists a constant matrix $L_0$, such that, $\eta_i$ will exponentially converge to $v$ as $t$ goes to infinity, in the sense of $||\eta_i-v||\leq c_0  e^{-\lambda_0 t}$ for some constants $c_0$ and $\lambda_0$.
\end{lemma}
{\bf Proof}:  Let $\eta=\text{col}(\eta_1,\dots,\eta_N)$ and $\bar \eta=\eta-\textbf{1}\otimes v$. After some mathematical manipulations, it follows
\begin{align*}
  \dot{\bar \eta}=(I_N\otimes S-H_{\sigma(t)}\otimes L_0F)\bar \eta
\end{align*}
Note that $H_{\sigma(t)}$ is constant and positive definite under Assumption \ref{ass:graph} during each interval. We first consider this problem in an interval $[t_i,t_{i+1})$. Assume $\sigma(t)=p$ for $t\in[t_i, t_{i+1})$, there exists a unitary matrix $U_p$ such that $\Lambda_p\triangleq U_p^TH_p U_p=\text{diag}\{\lambda_1^p,\dots,\lambda_N^p\}$. Let $\hat \eta=(U_p^T\otimes I_N)\bar \eta$, then,
\begin{align*}
  \dot{\hat \eta}=(I_N\otimes S-\Lambda_p\otimes L_0F)\hat \eta
\end{align*}
that is, $ \dot{\hat \eta}_i=(S-\lambda_i^p L_0F)\hat \eta_i$ for $i=1,\dots,N$. By selecting $L_0$ as that in Lemma \ref{lem:lmi} and letting $V_{\eta_i}={\hat \eta}_i^T P {\hat \eta}_i$, we can derive
$\dot{V}_{\eta_i}\leq -c V_{\eta_i}$. Recalling the dwell-time assumption, this inequality holds for all $t$. Let $V_{\eta}(t)=\sum_{i=1}^NV_{\eta_i}(t)$, it follows $\dot{V}_{\eta}\leq-c V_{\eta}$. Note that $\bar \eta^T \bar \eta=\hat \eta^T \hat\eta$, then
\begin{align*}
  ||\eta_i-v||^2\leq \bar \eta^T \bar \eta \leq \lambda_{min}(P)^{-1}{V}_{\eta}(t)\leq\lambda_{min}(P)^{-1}V_{\eta}(0)e^{-ct}.
\end{align*}
The conclusion follows readily. \hfill\rule{4pt}{8pt}

\begin{remark}
  It has been proved that, a sufficient and necessary condition for the linear matrix inequality \eqref{eq:lyapunov} is the detectability of $(F, S)$ (\cite{boyd1994linear}).  Hence, the observability of the leader system is sufficient to build an exponentially convergent distributed observer of the form \eqref{sys:observer} under Assumption \ref{ass:graph}.
\end{remark}
\begin{remark}
A similar result has been obtained in \cite{su2012tac} but under a fixed topology. Extensions under switching topologies were also obtained in \cite{su2012smc-cooperative} but assuming $S$ has no eigenvalues with positive real parts. Lemma \ref{lem:observer} allows that $S$ can have eigenvalues with positive real parts and hence includes unbounded references for this multi-agent system under switching topologies. Moreover, while full-information of the leader was used in \cite{su2012tac, su2012smc-cooperative}, i.e., $F=I_q$, we only needs the referenced output $y_0$ here and avoid the cases when the availability of the leader's full-information is not desirable.
\end{remark}

Next, we first consider a dynamic state feedback control to achieve the robust consensus tracking goal and then propose an output feedback design.

Let $\alpha(s)=s^{p_m}+\alpha_1 s^{p_m-1}+\dots+\alpha_{p_m-1}s+\alpha_{p_m}$ be the minimal polynomial of $S$, $G_1=I_{q}\otimes \beta$, and $G_2={I}_q\otimes \gamma$ with
\begin{align*}
\beta=\left[\begin{array}{c|c}
  0&I_{p_m-1}\\ \hline
  -\alpha_{p_m}&-\alpha_{p_m-1}\dots \alpha_{1}
\end{array}\right],\quad\gamma=\text{col}(\underbrace{0,0,\dots,0}_{p_m-1},1).
\end{align*}

Evidently, the pair $(G_1, G_2)$ is controllable. It is the so-called $q$-copy of $S$ (\cite{huang2004}).  By incorporating its $q$-copy of the leader, we propose a state feedback control law of the form
\begin{align}\label{ctr:state}
\begin{split}
u_i&=K_{1i}x_i+K_{2i}\xi_i\\
\dot{\xi}_i&=G_1\xi_i+G_2(y_i+F \eta_i)\\
\dot{\eta}_i&=S\eta_i-L_0F\eta_{vi}
\end{split}
\end{align}
where $G_1$, $G_2$ are as defined and $K_{1i}$, $K_{2i}$ to be determined later.

By a similar proof of Lemma 1.23 in \cite{huang2004}, under Assumption \ref{ass:trans}, the pair
\begin{align}\label{eq:ctr}
\left(\begin{bmatrix}
  A_i&\mathbf{0}_{n_i\times q p_m}\\
  G_2C_i&G_1
\end{bmatrix}, \begin{bmatrix}
  B_i\\
  G_2D_i
\end{bmatrix}\right)
\end{align}
is controllable for $i=1,\dots,N$. Hence, there exist constant matrices $K_{1i}$ and $K_{2i}$ with proper dimensions such that the matrix
\begin{align}\label{eq:closed-matrix}
A_{ci}\triangleq
\begin{bmatrix}
  A_i+B_iK_{1i}&B_iK_{2i}\\
  G_2(C_i+D_iK_{1i})&G_1+G_2D_iK_{2i}
\end{bmatrix}
\end{align}
is Hurwitz.

The regulator equations play a key role in the output regulation problem (\cite{huang2004}), and the following lemma establishes its solvability for the coordination design. Its proof is similar with that of Lemma 1.27 in \cite{huang2004} and thus omitted.
\begin{lemma}\label{lem:re}
  There exists a neighborhood $W_i$ of $w_i=0$, such that for any $w_i\in W_i$, the following equations:
  \begin{align}\label{eq:re}
    X_i(w)S&=A_i(w)X_i(w)+B_i(w)U_i(w),\\
    Z_i(w)S&=G_1Z_i(w)
  \end{align}
  have a solution $(X_i(w),\, Z_i(w))$ with $U_i(w)=K_{1i}X_i(w)+K_{2i}Z_i(w)$. Moreover, $X_i(w)$ and $Z_i(w)$ satisfy $$C_i(w)X_i(w)+ D_i(w)U_i(w)+F=0.$$
\end{lemma}

It is time to give our first main result.
\begin{theorem}\label{thm:state}
Under Assumptions \ref{ass:graph} and \ref{ass:trans}, the robust consensus tracking problem of this multi-agent system can be solved by the control law \eqref{ctr:state} with selected matrices $K_{1i}$, $K_{2i}$, $i=1,\dots,N$, and $L_0$.
\end{theorem}
{\bf Proof}: Under Assumption \ref{ass:trans} and by Lemma \ref{lem:re}, for any $i$,  there exists a neighborhood $W_i$ of $w_i=0$, such that for any $w_i\in W_i$, there exist two matrices $X_i(w)$ and $Z_i(w)$ satisfying
  \begin{align*}
    X_i(w)S&=A_i(w)X_i(w)+B_i(w)U_i(w)\\
    Z_i(w)S&=G_1Z_i(w)
  \end{align*}
where $U_i(w)=K_{1i}X_i(w)+K_{2i}Z_i(w)$.  Taking $W$ as the direct product of all $W_i$ (i.e., $W=W_1\times\dots \times W_N$) and performing a transformation $\bar x_i=x_i-X_i(w)v,\; \bar \xi_i=\xi_i-Z_i(w)v$ gives
\begin{align*}
  \dot{\bar{x}}_i&=(A_i(w)+B_i(w)K_{1i})\bar x_i+B_i(w)K_{2i}\bar \xi_i\\
  \dot{\bar{\xi}}_i&=G_2(C_i(w)+D_i(w)K_{1i})\bar x_i+(G_1+G_2D_i(w)K_{2i})\bar \xi_i-G_2\bar \eta_i\\
  e_i&=(C_i(w)+D_i(w)K_{1i})\bar x_i+D_i(w)K_{2i}\bar\xi_i
\end{align*}
Denote $\bar x_{ci}=\text{col}(\bar x_i,\,\bar\xi_i)$ and rewrite these equations in a compact form:
\begin{align*}
  \dot{\bar x}_{ci}&=A_{ci}(w)\bar x_{ci}+E_{ci}(w)\bar \eta_i\\
  e_i&=C_{ci}(w)\bar x_{ci}
\end{align*}
where $E_{ci}(w)=\text{col}(\textbf{0}_{n_i\times n_z},\,G_2F),\; C_{ci}(w)=[C_i(w)+D_i(w)K_{1i}, \, D_i(w)K_{2i}]$ are constant matrices with uncertain parameters.  By the exponential stability of $A_{ci}(w)$ and Lemma \ref{lem:observer}, we apply Lemmas 4.6 and 4.7 in \cite{khalil2001} and obtain the convergence of $\bar x_{ci}$. As a result, the proof is completed. \hfill\rule{4pt}{8pt}

\begin{remark}
  It can be found that if $y_0=v$~(i.e., $F=-I_q$), it means that the state variable $v$ of the leader can be directly obtained when it is reachable from some other agent. This special case has been partly considered in \cite{su2012tac, su2012smc-cooperative}. We extend the control law in \cite{su2012tac, su2012smc-cooperative} by only using the measurement output of the leader. Similar control laws were proposed in \cite{hong2006tracking,hong2009multi} when the followers are all integrators, while here we consider more general cases in this formulation.
\end{remark}

In some circumstances, the follower may not have access to its own state and only its output measurement is available. The following theorem shows how the robust consensus tracking problem can be solved by an output feedback control law.

\begin{theorem}\label{thm:output}
Under Assumptions \ref{ass:graph} and \ref{ass:trans}, there exist proper matrices $K_{1i}$, $K_{2i}$, $K_{3i}$, and $L_0$, such that the robust consensus tracking problem of this multi-agent system is solved by the following control law:
\begin{align}\label{ctr:output}
\begin{split}
u_i&=K_{1i}\zeta_i+K_{2i}\xi_i\\
\dot{\zeta}_i&=A_i\zeta_i+B_iu_i+K_{3i}(y_i-C_i\zeta_i-D_iu_i)\\
\dot{\xi}_i&=G_1\xi_i+G_2(y_i+F \eta_i)\\
\dot{\eta}_i&=S\eta_i-L_0F\eta_{vi}
\end{split}
\end{align}
\end{theorem}

{\bf Proof}: The proof is similar with that of Theorem \ref{thm:state}. Since the pair \eqref{eq:ctr} is controllable, we still take $K_{1i}$ and $K_{2i}$ as Theorem \ref{thm:state}. Note that $(C_i,\, A_i)$ is observable, there exists an $K_{3i}$ such that $A_i-K_{3i}C_i$ is Hurwitz. Denote $\hat x_{ci}=\text{col}(x_i,\,\zeta_i,\,\xi_i)$ and consider the closed-loop system of agent $i$ in a compact form as
\begin{align*}
  \dot{\hat {x}}_{ci}=\bar A_{ci}(w)\hat x_{ci}+\bar E_{ci}(w)\eta_i
\end{align*}
where $\bar E_{ci}(w)=\text{col}(\mathbf{0}_{n_i\times m}, \,\mathbf{0}_{n_i\times m},\,G_2F)$, and \begin{align*}
  \bar A_{ci}(w)=\begin{bmatrix}
A_i(w)&B_i(w)K_{1i}&B_i(w)K_{2i}\\
K_{3i}C_i(w)&\Xi_{1i}&\Xi_{2i}\\
G_2C_i(w)&G_2D_i(w)K_{1i}&G_1+G_2D_i(w)K_{2i}
  \end{bmatrix},
\end{align*}
with $\Xi_{1i}\triangleq A_i-K_{3i}C_i(w)+B_iK_{1i}+K_{3i}(D_i(w)-D)K_{1i}$ and $\Xi_{2i}\triangleq B_iK_{2i}+K_{3i}(D_i(w)-D)K_{2i}$.

The nominal system matrix $\bar A_{ci}(0)$ is
\begin{align*}
  \bar A_{ci}(0)=\begin{bmatrix}
A_i&B_iK_{1i}&B_iK_{2i}\\
K_{3i}C_i&A_i-K_{3i}C_i+B_iK_{1i}&B_iK_{2i}\\
G_2C_i&G_2D_iK_{1i}&G_1+G_2D_iK_{2i}
  \end{bmatrix},
\end{align*}
Let
\begin{align*}
  T_i=\begin{bmatrix}
    I_{n_i}&0&0\\
-I_{n_i}&I_{n_i}&0\\
0&0&I_{n_z}
  \end{bmatrix}.
\end{align*}
Then,
\begin{align*}
  \hat A_{ci}(0)&\triangleq T_i\bar A_{ci}T_i^{-1}\\
  &=\begin{bmatrix}
A_i+B_iK_{1i}&B_iK_{1i}&B_iK_{2i}\\
0&A_i-K_{3i}C_i&0\\
G_2(C_i+D_iK_{1i})&G_2D_iK_{1i}&G_1+G_2D_iK_{2i}
  \end{bmatrix}
\end{align*}
Recalling the selection of $K_{1i}$,\,$K_{2i}$,\, $K_{3i}$,  $\hat A_{ci}(0)$ and hence $\bar A_{ci}$ is Hurwitz. There exists a neighborhood $W$ of $w=0$, such that $\bar A_{ci}(w)$ is also Hurwitz.  Following the same procedure as Lemma \ref{lem:re}, there exist unique matrices $\bar X_i(w)$, $\bar Y_i(w)$, and $\bar Z_i(w)$ satisfying
  \begin{align*}
    \bar X_i(w)S&=A_i(w)\bar X_i(w)+B_i(w)\bar U_i(w)\\
    \bar Y_i(w)S&=A_i(w)\bar X_i(w)+B_i(w)\bar U_i(w)+K_{3i}C_i(w)(\bar X_i(w)-\bar Y_i(w))\\
    \bar Z_i(w)S&=G_1\bar Z_i(w)+G_2[C_i(w)\bar X_i(w)+D_i(w)\bar U_i(w)+F]
  \end{align*}
where $\bar U_i(w)=K_{1i}\bar Y_i(w)+K_{2i}\bar Z_i(w)$. From the uniqueness, $\bar X_i(w)=\bar Y_i(w)$. By the similar techniques of Lemma 1.27 in \cite{huang2004}, we can also obtain $C_i(w)\bar X_i(w)+D_i(w)\bar U_i(w)+F=0$.

Performing a transformation $\bar x_i=x_i-\bar X_i(w)v, \, \bar \zeta_i=\zeta_i-\bar Y_i(w)v,\,\bar \xi_i=\xi_i-\bar Z_i(w)v$ and denoting $\bar x_{ci}=\text{col}(\bar x_i,\,\bar \zeta_i,\,\bar \xi_i)$ gives
\begin{align*}
\dot{\bar{x}}_{ci}&=\bar A_{ci}(w)\bar x_{ci}+\bar E_{ci}(w)\bar \eta_i\\
e_i&=C_{ci}(w)\bar x_{ci}
\end{align*}
where $C_{ci}(w)=[C_i(w),\,D_i(w)K_{1i}, \, D_i(w)K_{2i}]$ are constant matrices with uncertainties.  By the stability of $\bar A_{ci}(w)$ and Lemma \ref{lem:observer}, we apply Lemmas 4.6 and 4.7 in \cite{khalil2001} and obtain the convergence of $\bar x_{ci}$. As a result, the proof is completed.   \hfill\rule{4pt}{8pt}

\begin{remark}\label{rem:su-observer}
  This problem has been partially investigated in \cite{su2012smc-cooperative} under switching topologies. However, the constructed controllers heavily relied on the exact known system matrices and hence can't admit uncertainties. The internal model-based control law here facilitates us to allow small parameter uncertainties and thus obtain a robust control law. Moreover, without restricting the eigenvalues of $S$ are in the closed left half-plane, the common Lyapunov function technique helps us handle the switching topologies even the references may be unbounded.
\end{remark}
\begin{remark}\label{rem:wang-canonical}
    Another relevant paper is \cite{wang2011robust}, in which a robust output regulation problem of linear multi-agent systems was studied by a so-called canonical internal model. The main differences between our work and that of \cite{wang2011robust} are at least two-fold. First, the method proposed in \cite{wang2011robust} only applies to the case when the outputs have the same dimension with its inputs. Second, only an output feedback control law was derived in \cite{wang2011robust}, while here both state feedback and output feedback control laws are proposed. Moreover, since the out feedback control law in \cite{wang2011robust} needs to estimate the state of the internal model, its dimension may be much higher than the control \eqref{ctr:output} when $q>1$.
\end{remark}

\section{Simulations}

\begin{figure}
  \centering
  \subfigure[The graph $\mathcal{G}_1$]
    {\centering
\begin{tikzpicture}[shorten >=1pt, node distance=1.2 cm, >=stealth',
every state/.style ={circle, minimum width=0.2cm, minimum height=0.2cm}, auto]
\node[align=center,state](node0) {0};
\node[align=center,state](node1)[right of=node0]{1};
\node[align=center,state](node2)[right of=node1]{2};
\node[align=center,state](node3)[right of=node2]{3};
\path[->]   (node0) edge (node1)
            (node1) edge [bend right] (node2)
            (node2) edge [bend right] (node1)
            (node0) edge [bend left]  (node3)
            ;
\end{tikzpicture}
  }\quad
  \subfigure[The graph $\mathcal{G}_2$]
   {\centering
\begin{tikzpicture}[shorten >=1pt, node distance=1.2 cm, >=stealth',
every state/.style ={circle, minimum width=0.2cm, minimum height=0.2cm}, auto]
\node[align=center,state](node0) {0};
\node[align=center,state](node1)[right of=node0]{1};
\node[align=center,state](node2)[right of=node1]{2};
\node[align=center,state](node3)[right of=node2]{3};
\path[->]   (node0) edge (node1)
            (node0) edge [bend left](node2)
            (node2) edge [bend right] (node3)
            (node3) edge [bend right] (node2)
            ;
\end{tikzpicture}
    }
\caption{The communication graphs.}\label{fig:graph}
\end{figure}
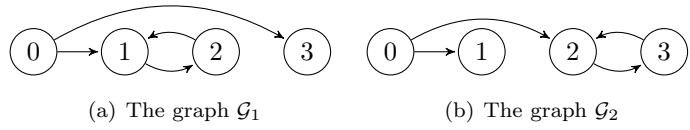

For illustrations, we present an example and consider three followers with following system matrices
\begin{align*}
  A_1=1+\epsilon_1,\; B_1=1+\epsilon_1,\; C_1=1+\epsilon_1,\; D_1=1+\epsilon_1,
\end{align*}
\begin{align*}
  A_2=\begin{bmatrix}
  0&1+\epsilon_2 \\-1+\epsilon_2&0
\end{bmatrix},\; B_2=\mbox{col}(0,1+\epsilon_2),\; C_2=[1+\epsilon_2,\;0],\; D_2=\epsilon_2,
\end{align*}
\begin{align*}
  A_3=\begin{bmatrix}
  \epsilon_3&1&0\\
  -1+\epsilon_3&0&1\\
  2&\epsilon_3&1\\
\end{bmatrix},\; B_3=\mbox{col}(0,1,1+\epsilon_1),\; C_3=[0,\;1,\;0],\; D_3=1.
\end{align*}
The leader is described by
\begin{align*}
  \dot{v}_1=v_2,\; \dot{v}_2=-v_1,\; y_0=-v_1.
\end{align*}
We assume the interconnection topology is switching between graph $\mathcal{G}_i$ $(i=1,2)$ described by Figure \ref{fig:graph}. The switchings are periodically carried out as $\{\mathcal{G}_1, \mathcal{G}_2, \mathcal{G}_1, \mathcal{G}_2, \cdots\}$ with periods $t=5$.

\begin{figure}
  \centering
  \subfigure[The trajectories of the leader and followers. ]{
    \includegraphics[width=0.42\textwidth]{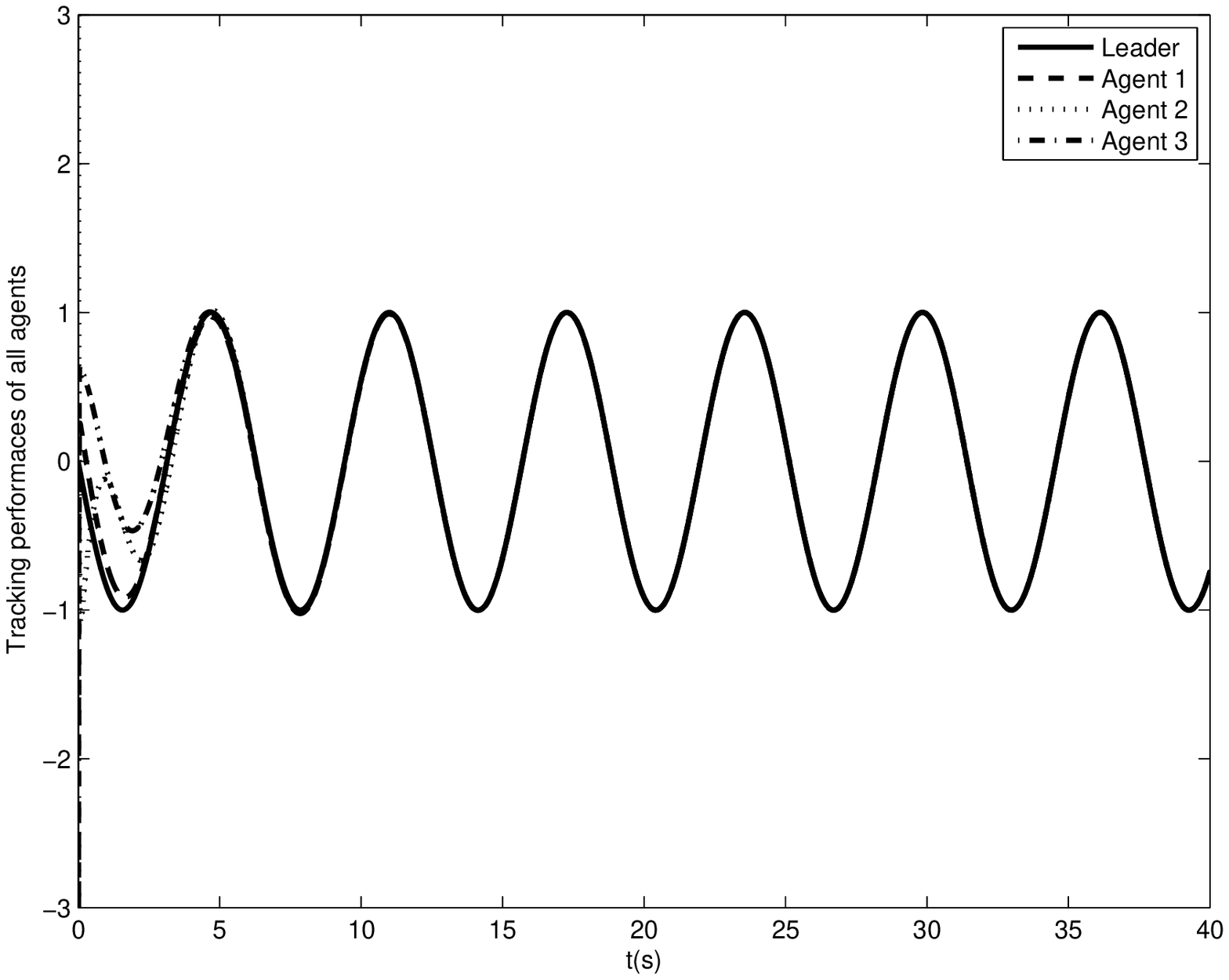}
  }
  \subfigure[The tracking errors of three followers.]{
    \includegraphics[width=0.42\textwidth]{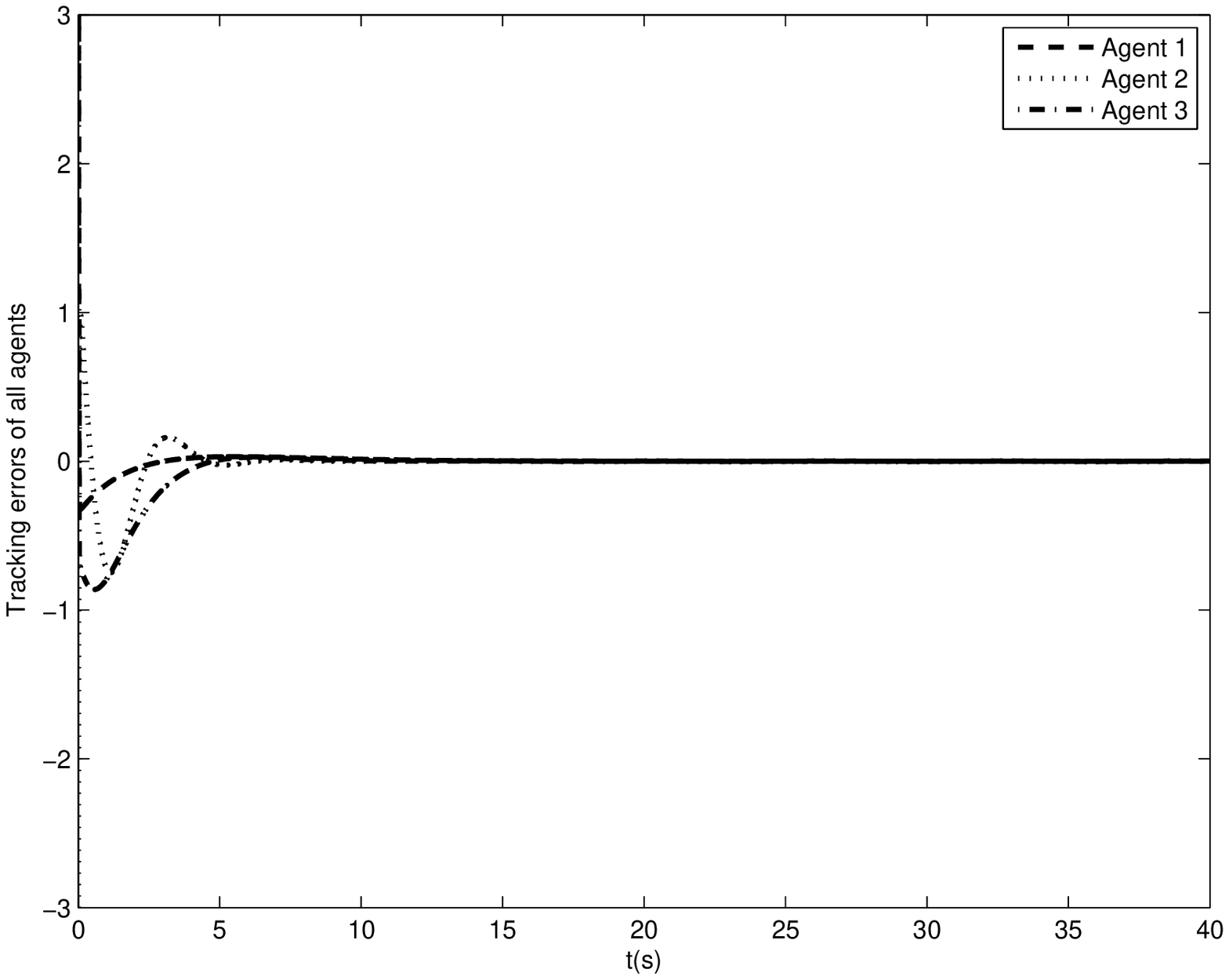}
    }
\caption{The tracking performance of three agents under the controller \eqref{ctr:state}. }\label{fig:simu1}
\end{figure}

\begin{figure}
  \centering
  \subfigure[The trajectories of the leader and followers.]{
    \includegraphics[width=0.42\textwidth]{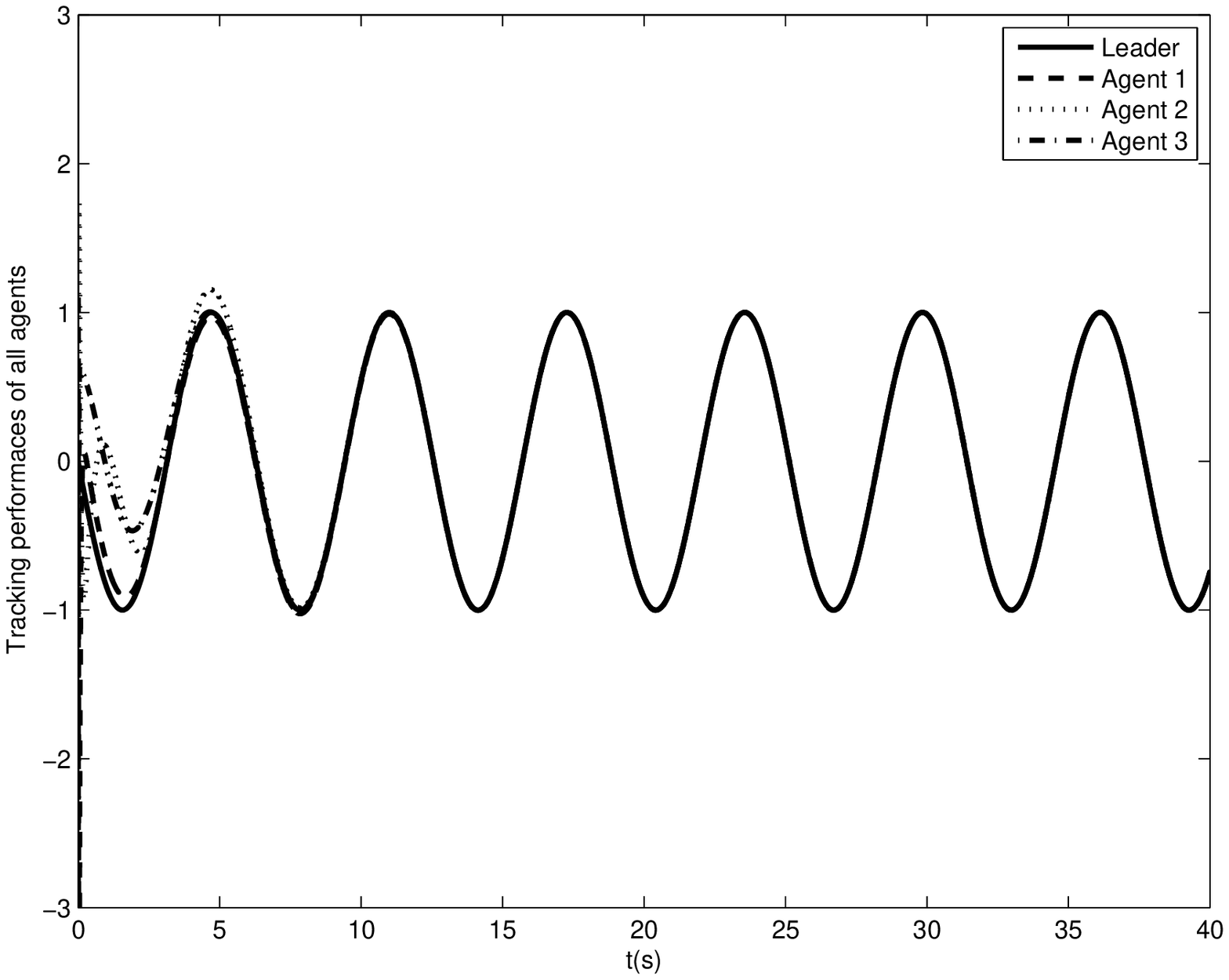}
  }
  \subfigure[The tracking errors of three followers.]{
    \includegraphics[width=0.42\textwidth]{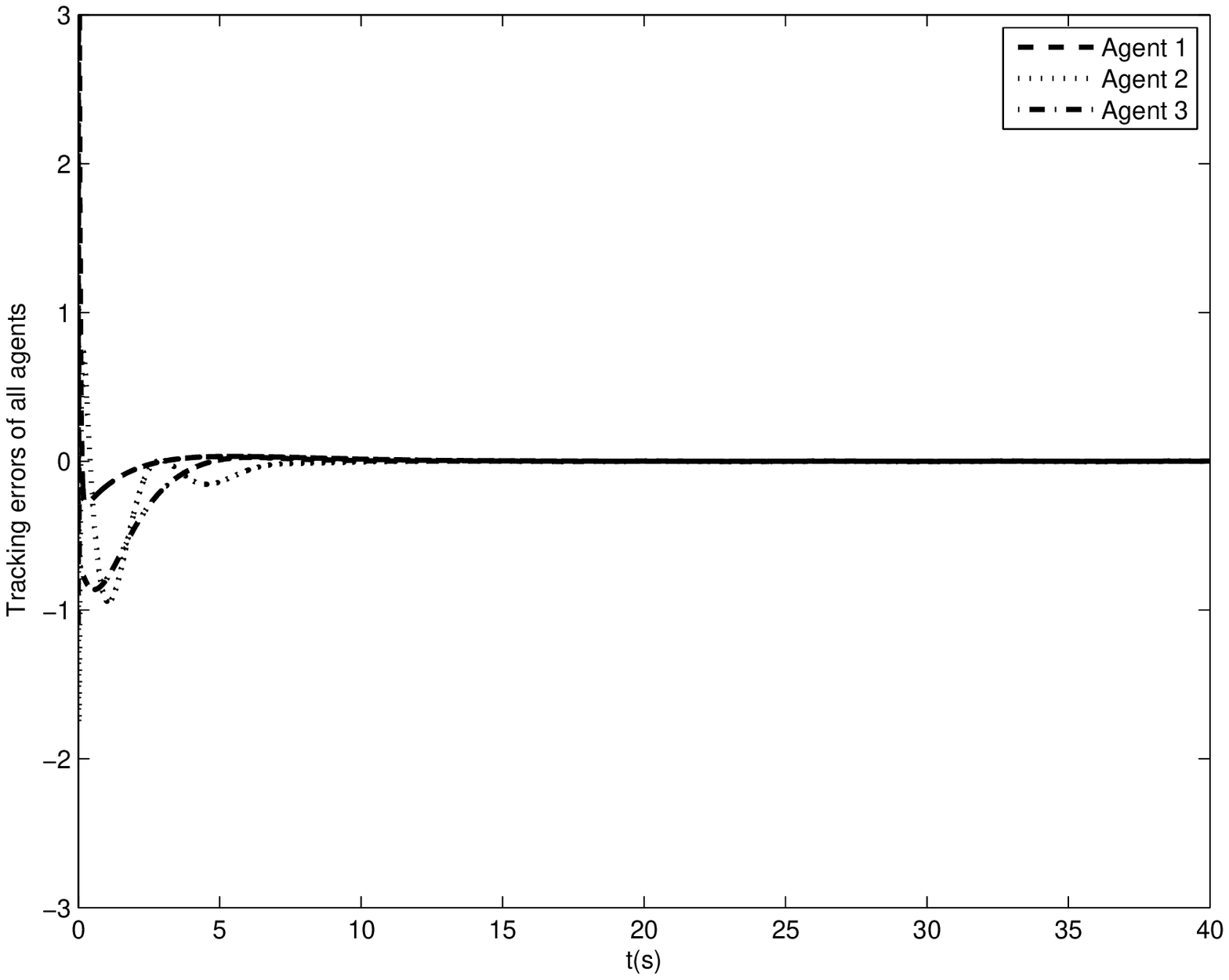}
    }
\caption{The tracking performance under the controller \eqref{ctr:output}.}\label{fig:simu2}
\end{figure}

The uncertain parameters $\epsilon_1, \epsilon_2, \epsilon_3$ are taken between $[-1,1]$. The robust consensus tracking performances under state and output feedback control laws were depicted in Figures \ref{fig:simu1} and \ref{fig:simu2}.

\section{Conclusions}
A robust consensus tracking problem of heterogeneous multi-agent systems was solved to admit small uncertain parameters in the agents' systems. In conjunction with the internal model techniques, a common Lyapunov function was used to overcome the challenges of switching topologies.  Future work includes extensions to more general graphs and multiple leaders.

\end{document}